\documentclass[a4paper,10pt,epsfig,graphicx,subfigure,pstricks,twoside]{article}
\usepackage[dvips]{graphics,color}
\usepackage[english]{babel}
\usepackage{epsf}
\usepackage{verbatim}
\usepackage{stmaryrd}
\usepackage{latexsym}
\usepackage{makeidx}
\usepackage{epsfig}
\usepackage{amsfonts}
\usepackage{amsmath}
\usepackage{amssymb}

\input xy
\xyoption{all}
\makeindex
\newtheorem{Def}{Definition}[section]
\newtheorem{Prop}[Def]{Proposition}

\newtheorem{Lem}[Def]{Lemma}
\def\norm#1{\left |#1\right |}
\def\Norm#1{\left \|#1\right \|}

\def\endproof{\ \hfill\hbox{\vbox{\hrule\hbox{\vrule
height5pt\kern5pt\vrule height5pt}\hrule}}\par\medskip\rm}

\title{\bf Convergence to the time average \\
by stochastic regularization}
\author{ {O. Bernardi, F. Cardin, M. Guzzo }\\ \\
Dipartimento di Matematica \\
Universit\`a degli Studi di Padova \\
Via Trieste, 63 - 35121 Padova, Italy \\}
\date{ \ }
\begin{document}
\maketitle
\begin{abstract} \par \noindent  
We compare the rate of convergence to the time average of a function over an
integrable Hamiltonian flow with the one obtained by a stochastic
perturbation of the same flow. Precisely, we provide detailed estimates in
different Fourier norms and we prove the convergence even in a Sobolev norm
for a special vanishing limit of the stochastic perturbation.
\end{abstract}

\section{Introduction}\label{intro}

The time averages of functions with respect to the flow of Hamiltonian 
systems are extensively studied in Ergodic Theory and Hamiltonian Perturbation
Theory. In particular, the averages over integrable flows 
are commonly used as generating 
functions of averaging canonical transformations. In this setting 
it is well known since Poincar\'e that resonances related to the so--called 
small divisors 
represent topological obstructions to the regularity of the time averages in 
open sets of the phase--space. The celebrated KAM and Nekhoroshev 
Theorems (\cite{Kolmogorov}, \cite{Arnolda}, \cite{Moser}, \cite{Nekh77})
overcome this problem with a 
refined use of algebraic as well as geometric treatment of small divisors. 
More recently, the so--called weak KAM theories (see \cite{FATHI}, 
\cite{MATHER}, \cite{EV09}) have studied the problem by new perspectives, 
based on variational and PDE regularizations by viscosity techniques.  
\\
\indent In this paper, in order to estimate the rate of convergence to the
time 
average, we exploit a correspondence between standard viscosity 
regularizations of PDEs (see for example \cite{FW}) and the averages of 
functions over stochastic perturbations. We prove that a vanishing  
stochastic regularization of the time average over an integrable flow 
converges to the time average in a Sobolev norm accounting the 
first derivatives. Precisely,
let us consider the integrable Hamiltonian system with Hamilton 
function $H(I,\varphi):=h(I)$, defined on the action--angle 
phase--space $A\times {\Bbb T}^n$, where  $A\subseteq  \mathbb{R}^n$ is 
open bounded and $g(I) := \nabla h(I)$ is a diffeomorphism over $A$
such that  
\begin{equation} 
\norm{g(I)} \le C,\qquad \max_{i,j}\norm{{\partial g_i \over \partial I_j}(I)} 
\leq D, \qquad \left| \text{det} 
\frac{\partial g}{\partial I}(I) \right| \geq m 
 \label{GI}
\end{equation}
$\forall I \in A$,  for some positive 
constants $C,D,m> 0$. We will also denote by $\lambda>0$ a Lipschitz constant 
for $g$ in the set $A$. 
 \\
\indent For any smooth phase--space function $f(I,\varphi)$, we consider its 
 finite--time average
\begin{equation} \label{GT}
G^T(I,\varphi) := \frac{1}{T}\int^T_0 f(\phi^{t}(I,\varphi)) dt  ,
\end{equation}
where $\phi^t(I,\varphi) = (I,\varphi + g(I) t)$. By denoting with 
$$
f(I,\varphi) := \sum_{k\in {\Bbb Z}^n} f_k(I)e^{ik\cdot \varphi}, \qquad
G^T(I,\varphi) := \sum_{k\in {\Bbb Z}^n} G^T_k(I)e^{ik\cdot \varphi}
$$
the Fourier expansions of $f$ and $G^T$, we have 
\begin{equation*} \label{intro1}
G^T_k(I) = \begin{cases} f_k(I) \qquad & \text{if } k \cdot g(I) = 0 \\ \\
\displaystyle{f_k(I) \frac{e^{ik \cdot g(I) T} - 1}{i k \cdot g(I) T}} \qquad & \text{if } k \cdot g(I) \ne 0
\end{cases}
\end{equation*}
With evidence, if $f_k \ne 0$ for a suitably large --that is generic-- set of indices 
$k\in {\Bbb Z}^n$, the presence of small divisors $k \cdot g(I)$ represents an obstruction to the regularity both for $G^T$ and for its 
limit
\begin{equation} \label{limite}
\bar{f}(I,\varphi) := \lim_{T \to +\infty} G^T(I,\varphi)  .
\end{equation}
We remark that  the Fourier coefficients  $G^T_k(I)$ are similar to 
 the Fourier coefficients  of
\begin{equation*}
\chi(I,\varphi)= - \sum_{k \in {\Bbb Z}^n\backslash 0} {f_k(I)\over 
i k\cdot g(I)}e^{ik \cdot \varphi},
\end{equation*}
whose $\epsilon$--time flow $\phi^{\epsilon}_{\chi}$ formally conjugates the 
quasi--integrable Hamiltonian system 
\begin{equation*}
H_\epsilon(I,\varphi) = h(I)+\epsilon f(I,\varphi)
\end{equation*}
to its first order average
$$
(H_\epsilon\circ \phi^{\epsilon}_{\chi}) (I,\varphi) =h(I)+\epsilon f_0(I)+
{\cal O}(\epsilon^2) .
$$
Of course, $\chi$ and $G^T$ are affected by the same convergence problems. \\
\indent We assume from now on that $f$ is smooth and with generic Fourier 
expansion. Precisely, let us introduce for any $k\in {\Bbb Z}^n$ the 
resonant manifold
\begin{equation}
{\cal R}_k = \{ I\in A:\ \ k\cdot g(I)=0\}  ,
\label{resman}
\end{equation}
as well as
\begin{equation}
{\cal R}_k(f) = \{ I\in A:\ \ k\cdot g(I)=0\ \ {\rm and}\ \ \norm{f_k(I)}>0
\}  .
\label{resmanif}
\end{equation}
Then, we assume that the set
\begin{equation}
{\cal R}(f) = \bigcup_{k \in {\Bbb Z}^n\backslash 0}{\cal R}_k(f)
\label{erre}
\end{equation}
is dense in $A$.\\ 
\indent We now consider the regularization of $G^T$ based on a vanishing 
stochastic perturbation, previously introduced in \cite{BCGZ} by following 
a technique described in \cite{FW}. More 
precisely, let $(\Omega, \mathcal{F}, P)$ be a probability space and 
$w_t: \Omega \to \mathbb{R}^n$ a $n$--dimensional Wiener process. Then, we 
obtain a stochastic differential equation by perturbing the Hamilton equations 
with a white noise 
\begin{equation} \label{stochintro}
\begin{cases} \dot{I}_t = 0 \\ 
\dot{\varphi}_t = g(I) + 2\nu \dot{w}_t
\end{cases}
\end{equation}
whose flow is $\Phi^t_{\nu}(I,\varphi,\omega) = (I,\varphi +
 g(I)t + 2\nu w_t(\omega))$. As in \cite{BCGZ}, for $\mu,\nu > 0$ we introduce 
\begin{equation} \label{me numu}
F^{\mu,\nu}(I,\varphi) := \mu M_{(I,\varphi)}\left( \int^{+\infty}_0 f(\Phi^t_{\nu}(I,\varphi,\omega))e^{-\mu t} dt \right).
\end{equation}
In the previous formula, $M_{(I,\varphi)}$ represents, for $(I,\varphi)$ fixed, the average on
all the trajectories of the Brownian motion (\ref{stochintro}), while the 
exponential damping $e^{-\mu t}$  allows us to
interpret $F^{\mu,\nu}$ as an effective average over 
a time interval of some multiples of $1/\mu$ (see \cite{BCGZ}). Moreover, 
in this paper this factor will play an essential role to ensure 
the convergence for $(\mu,\nu)\rightarrow (0,0)$.  \\ 
\indent In order to study the convergence properties of $G^T$ and $F^{\mu,\nu}$ to the time average $\bar{f}$, we introduce specific norms on 
$A\times {\Bbb T}^n$. In more detail, 
for any function $u(I,\varphi) = \sum_{k \in \mathbb{Z}^n} u_k(I) e^{ik \cdot \varphi}$ on $A\times {\Bbb T}^n$, the 
uniform Fourier norm 
\begin{equation} \label{norm1-intro}
|u|^{\infty} := \sum_{k \in \mathbb{Z}^n} \sup_{I \in A} |u_k(I)|
\end{equation}
as well as the norms obtained with averages over the action space 
\begin{eqnarray} \label{norm2-intro}
|u|^0 := \sum_{k \in \mathbb{Z}^n} \int_A |u_k(I)| dI
\end{eqnarray}
and
\begin{eqnarray} \label{norm3-intro}
|u|^1 := |u|^0 +\sum_{k \in \mathbb{Z}^n}
\sum_{j=1}^n \int_{A} \left (
\left| \frac{\partial u_k}{\partial I_j}(I) \right|
+\left| k_j u_k(I)\right| \right ) dI  .
\end{eqnarray}
Let us remark that, by considering the usual $L^1$ and Sobolev $W^{1,1}$ norms 
on $A\times {\Bbb T}^n$, in particular
$$
\Norm{u}_{W^{1,1}}= \Norm{u}_{L^1}+
\sum_{j=1}^n 
\left (
\Norm{ \frac{\partial u}{\partial I_j}}_{L^1}+
\Norm{ \frac{\partial u}{\partial \varphi_j}}_{L^1}\right ) 
$$
and we have
$$
{1\over (2\pi)^n}
\Norm{u}_{W^{1,1}} \leq  \norm{u}^1 \leq 
{1\over (2\pi)^n}
\sum_{k \in \mathbb{Z}^n}
\Norm{u_k(I) e^{ik\cdot \varphi}}_{W^{1,1}}  
$$
\indent The paper presents the following results. In Proposition \ref{PROP2}, we first 
prove that for a generic $f$, both $G^T$ and $F^{\mu,\nu}$ do not 
converge to $\bar f$ in the uniform norm $\norm{\ \cdot }^\infty$, but they 
converge to 
$\bar f$ in the $\norm{\ \cdot \ }^0$ norm. The main 
result, given in Proposition \ref{PROP1}, concerns with the stronger norm 
 $\norm{\ \cdot \ }^1$: while the finite time average $G^T$ does not converge 
to $\bar f$ in the $\norm{\ \cdot \ }^1$ norm, we have 
$$
\lim_{i \rightarrow +\infty } |F^{\mu_i, \nu_i} - \bar{f}|^1 = 0
$$
for any sequence $\mu_i, \nu_i$ converging to zero and  such that $\lim_{i \to +\infty} \frac{\mu_i}{\nu_i} = 0$. \\
\indent The paper is organized as follows. In Section \ref{SE2} we state Propositions 
\ref{PROP2} and \ref{PROP1} on the convergence of the regularized
averages. Section \ref{SE3} is devoted to proofs.

\section{Convergence results} \label{SE2}

Let us consider a phase--space function $f(I,\varphi)$, its finite time 
average $G^T$ and its time average $\bar f$ defined in 
(\ref{GT}) and (\ref{limite}) respectively. In \cite{BCGZ} we have introduced two different approximations of $G^T$ offering 
a better notion of approximated first integral. The first one is
\begin{equation} \label{me nu}
F^{\mu}(I,\varphi) := \mu \int^{+\infty}_0 f(\phi^{t}(I,\varphi)) e^{-\mu t} dt  ,
\end{equation}
with $\mu=1/T$, while the second one is $F^{\mu,\nu}$, whose definition 
recalled in Section 
\ref{intro} (see (\ref{me numu})), comes from the above stochastic setting. 
Let us remark that 
$F^\mu$ represents an intermediate step between $G^T$ and $F^{\mu,\nu}$, in the 
sense that it is an exponentially damped average of $f$ with respect to 
the integrable flow, that is, $F^\mu = F^{\mu,0}$. However, the  improvement 
in 
the convergence properties to $\bar f$ is obtained only for $\nu>0$ (see the 
propositions below).\\ 
\indent We first discuss the convergence of the above 
approximated first integrals $G^T$, $F^{\mu}$ and $F^{\mu, \nu}$ to the time
average $\bar{f}$ both in the uniform Fourier norm $\norm{\ \cdot\ }^\infty$
 --see (\ref{norm1-intro})-- and in the action--averages based norm 
$\norm{\ \cdot\ }^0$ given in (\ref{norm2-intro}). In particular, we prove the 
next
\begin{Prop} \label{PROP2} 
The functions $G^T$, $F^{\mu}$ and
  $F^{\mu, \nu}$ do not uniformly Fourier converge to $\bar{f}$ in any set 
$B\times {\Bbb T}^n$ with $B\subseteq A$ open, 
but converge  to $\bar{f}$ in the $\norm{\ \cdot\ }^0$ norm on $A\times {\Bbb
  T}^n$. Precisely, we have
\begin{equation}
|G^T - \bar{f}|^0 \le {4  C^{n-1} \over m} 
\sum_{k\in {\Bbb Z}^n\backslash 0}|f_k|^{\infty}\  {3 + \log(\|k\| T C)\over T}
\label{convgt}
\end{equation}
\begin{equation}
|F^{\mu} - \bar{f}|^0 \le  {2  C^{n-1} \over m} 
\sum_{k\in {\Bbb Z}^n\backslash 0 }|f_k|^{\infty}\ \frac{\mu }{\|k\|} \left( 1 + \log {\|k\| C \over \mu }\right )
\label{convfm}
\end{equation}
\begin{equation}
 |F^{\mu,\nu} - \bar{f}|^0 \le  {2  C^{n-1} \over m} 
\sum_{k\in {\Bbb Z}^n\backslash 0 }|f_k|^{\infty}\ \frac{\mu }{\|k\|} \left( 1 + \log {\|k\| C \over \mu+\nu\Norm{k}^2 }\right ) .
\label{convfmn}
\end{equation}

\end{Prop} 
As it arises from the previous proposition, the three different finite--time approximations 
$G^T$, $F^{\mu}$ and $F^{\mu, \nu}$ behave in the same way with respect to 
the $\norm{\ \cdot\ }^\infty$ and $\norm{\ \cdot\ }^0$ norms. Indeed, the
difference consists in the convergence in 
the $\norm{\ \cdot\ }^1$ norm given in (\ref{norm3-intro}). In such a case, 
the $G^T$, $F^{\mu}$ do not converge to $\bar f$, and it is remarkable 
that the convergence of $F^{\mu,\nu}$ is obtained only in a special limit of 
 vanishing stochastic perturbation, as stated in the proposition below. 
\begin{Prop} \label{PROP1} The functions $G^T$ and $F^{\mu}$ 
do not converge to $\bar{f}$ in the $\norm{\ \cdot\ }^1$ norm
 on any set $B\times {\Bbb T}^n$ with $B\subseteq A$ open. 
Differently, for any $\mu,\nu>0$ the function $F^{\mu,\nu}$ satisfies
\begin{equation}
\norm{F^{\mu,\nu} -\bar{f}}^1 
\leq  {2 C^{n-1}\over m} 
\sum_{k\in {\Bbb Z}^n \backslash 0} \Bigg (  
\mu \Big [ 1 + \log {\|k\| C \over \mu+\nu\Norm{k}^2 }\Big ] 
\Big ( (1+n) \norm{f_k}^\infty  + \sum_{j=1}^n  \left|\frac{\partial
    f_k}{\partial I_j}\right|^{\infty} \Big )+
 {1\over 2} n^2 \pi D  {\mu \over \mu+\nu \Norm{k}^2} \norm{f_k}^\infty\Bigg )
\label{convfmnorm}
\end{equation}
on $A\times {\Bbb T}^n$. In particular, for any 
sequence $\mu_i$, $\nu_i>0$ converging to zero and such that 
$$\lim_{i \to +\infty} \frac{\mu_i}{\nu_i} = 0,$$ 
we have 
$$
\lim_{i \to + \infty} |F^{\mu_i, \nu_i} - \bar{f}|^1 = 0  .
$$
\end{Prop}
Let us remark that the convergence of $F^{\mu,\nu}$ to $\bar f$ requires 
a restriction of the sub-sequences $\mu_i,\nu_i$  because 
in (\ref{convfmnorm}) we find contributions proportional to 
${\mu \over \mu+\nu \Norm{k}^2}$, while the contributions $\mu 
\log {\|k\| C \over \mu+\nu\Norm{k}^2 }$ which are  dominant in
(\ref{convfmn}) converge for $(\mu,\nu) \to (0,0)$.

\indent The proofs of Propositions \ref{PROP2}, \ref{PROP1} are reported in 
Section \ref{SE3}.

\section{Proofs} \label{SE3}

The different time averages (\ref{GT}), (\ref{me nu}) and (\ref{me numu}) can be alternatively expressed in terms of their Fourier coefficients, as discussed in the following technical
\begin{Lem} \label{lemma} Let us consider 
\begin{equation} \label{effe}
f(I,\varphi) = \sum_{k \in \mathbb{Z}^n} f_k(I) e^{i k \cdot \varphi}.
\end{equation}
The Fourier coefficients of 
$$G^T(I,\varphi) = \sum_{k \in \mathbb{Z}^n} G^T_k(I) e^{i k \cdot \varphi}, \quad F^{\mu}(I,\varphi) = \sum_{k \in \mathbb{Z}^n} F^{\mu}_k(I) e^{i k \cdot \varphi}, \quad F^{\mu,\nu}(I,\varphi) = \sum_{k \in \mathbb{Z}^n} F^{\mu,\nu}_k(I) e^{i k \cdot \varphi}$$ 
are respectively
\begin{eqnarray} \label{coeff1}
G^T_k(I) = \begin{cases} f_k(I) \qquad & \text{if } k \cdot g(I) = 0 \\ \\
\displaystyle{f_k(I) \frac{e^{ik \cdot g(I) T} - 1}{i k \cdot g(I) T}} \qquad & \text{if } k \cdot g(I) \ne 0
\end{cases}
\end{eqnarray}
\begin{equation} \label{coeff2}
F^{\mu}_k(I) = -\mu \frac{f_k(I)}{ik \cdot g(I) - \mu}
\end{equation} 
and
\begin{equation} \label{coeff3}
F^{\mu,\nu}_k(I) = -\mu \frac{f_k(I)}{ik \cdot g(I) - \mu - \nu \|k\|^2}
\end{equation}
\end{Lem}
\textit{Proof.} The first equality easily follows from (\ref{GT}) and (\ref{effe}) by direct calculations. Indeed
\begin{eqnarray*}
G^T(I,\phi) &=& \frac{1}{T} \int^T_0 f(\phi^t(I,\varphi)) dt = \frac{1}{T} \int^T_0 \sum_{k \in \mathbb{Z}^n} f_k(I) e^{i k \cdot \varphi} e^{i k \cdot g(I) t} dt \\
&=& \frac{1}{T} \sum_{k \in \mathbb{Z}^n} f_k(I) e^{i k \cdot \varphi} \int^T_0 e^{i k \cdot g(I) t} dt.
\end{eqnarray*}
Moreover, from
$$\frac{1}{T}\int^T_0 e^{i k \cdot g(I) t} dt =
\begin{cases}
1 & \qquad \text{if } k \cdot g(I) = 0 \\ \\
\displaystyle{\frac{e^{i k \cdot g(I) T} - 1}{i k \cdot g(I) T}} & \qquad \text{if } k \cdot g(I) \ne 0 
\end{cases}$$
we immediately obtain formula (\ref{coeff1}). Similarly for (\ref{coeff2})
\begin{eqnarray*}
F^{\mu}(I,\varphi) &=& \mu \int^{+\infty}_0 f(\phi^t(I,\varphi)) e^{-\mu t}dt = \mu \int^{+\infty}_0 \sum_{k \in \mathbb{Z}^n}f_k(I) e^{ik \cdot \varphi}e^{ik \cdot g(I)t - \mu t}dt \\
&=& \mu \sum_{k \in \mathbb{Z}^n} f_k(I) e^{ik \cdot \varphi} \int^{+\infty}_0 e^{(ik \cdot g(I) - \mu)t} dt = -\mu \sum_{k \in \mathbb{Z}^n} \frac{f_k(I)}{ik \cdot g(I) - \mu} e^{ik \cdot \varphi}.
\end{eqnarray*}
We conclude by proving the equality (\ref{coeff3}). We first take into account 
(\ref{me numu}), so that
\begin{equation*} 
\int^{+\infty}_0 f(\Phi^t_{\nu}(I,\varphi,\omega))e^{-\mu t} dt = \sum_{k \in \mathbb{Z}^n} f_k(I) e^{i k \cdot \varphi} \int_0^{+\infty} e^{(i k \cdot g(I) - \mu) t} e^{2 i \nu k \cdot w_t(\omega)} dt.
\end{equation*}
As a consequence --see (\ref{me numu})-- we obtain 
\begin{eqnarray}
F^{\mu,\nu}(I,\varphi) &=& \mu M_{(I,\varphi)}\left( \int^{+\infty}_0 f(\Phi^t_{\nu}(I,\varphi,\omega))e^{-\mu t} dt \right) \cr
&=& \mu \sum_{k \in \mathbb{Z}^n} f_k(I) e^{i k \cdot \varphi} \int_{\Omega} \left[\int_0^{+\infty} e^{(i k \cdot g(I) - \mu) t} e^{2 i \nu k \cdot w_t(\omega)} dt\right] P(d\omega) \cr
&=& \mu \sum_{k \in \mathbb{Z}^n} f_k(I) e^{i k \cdot \varphi} \int_0^{+\infty} \left[ e^{(i k \cdot g(I) -\mu)t} \int_{\Omega} e^{2 i \nu k \cdot w_t(\omega)} P(d\omega) \right] dt.
\label{eqfmn}
\end{eqnarray}
Since $w_t: \Omega \to \mathbb{R}^n$ is a $n$--dimensional Wiener process, the corresponding covariance matrix $R(t) = R_{ij}(t) = t \delta_{ij}$ and therefore 
$$\int_{\Omega} e^{2 i \nu k \cdot w_t(\omega)} P(d\omega) = e^{-\nu \| k \|^2 t}.$$
Therefore, from equation (\ref{eqfmn}) we have
\begin{eqnarray*}
F^{\mu,\nu}(I,\varphi) &=& \mu \sum_{k \in \mathbb{Z}^n} f_k(I) e^{i k \cdot \varphi} \int^{+\infty}_0 e^{(ik \cdot g(I) - \mu -\nu \| k \|^2) t} dt = -\mu \sum_{k \in \mathbb{Z}^n} \frac{f_k(I)}{ik \cdot g(I) - \mu - \nu \|k\|^2} e^{ik \cdot \varphi}.
\end{eqnarray*}
\hfill $\Box$ \\
The next sections are devoted to the convergence results, in three different
norms, of $G^T$, $F^{\mu}$ and $F^{\mu, \nu}$ to the time average $\bar{f}$ 
defined in (\ref{limite}). From (\ref{coeff1}), we immediately obtain 
$\bar{f} = 
\sum_{k \in \mathbb{Z}^n} \bar{f}_k(I) e^{ik \cdot \varphi}$, with 
\begin{eqnarray} \label{coeff4}
\bar{f}_k(I) = \begin{cases} f_k(I) \qquad & \text{if } k \cdot g(I) = 0 \\ 
0 \qquad & \text{if } k \cdot g(I) \ne 0
\end{cases}
\end{eqnarray}

\subsection{Proof of Proposition \ref{PROP2}}

We start by proving that $G^T$ does not converge to $\bar{f}$ 
in the uniform Fourier norm.  Let us consider 
$$(G^T - \bar{f})(I,\varphi) := \sum_{k \in \mathbb{Z}^n} (G^T - \bar{f})_k(I)e^{i k \cdot \varphi}.$$
From (\ref{coeff1}) and (\ref{coeff4}) we immediately obtain 
\begin{equation} \label{differenza1}
(G^T - \bar{f})_k(I) = \begin{cases} 0 & \qquad \text{if } k \cdot g(I) = 0 \\ 
\displaystyle{f_k(I)\frac{e^{ik \cdot g(I) T} - 1}{i k \cdot g(I)T}} & \qquad \text{if } k \cdot g(I) \ne 0
\end{cases}
\end{equation}
Since the set ${\cal R}(f)$ defined in (\ref{erre}) is dense, there exists 
a dense set of points $ \bar I\in A$ such that  $\bar k \cdot g(\bar I)=0$
and $\norm{f_{\bar k}(\bar I)}>0$ for some $\bar k \in {\Bbb Z}^n\backslash
0$. Since $g$ is a diffeomorphism, we have
$$
\lim_{I \notin {\cal R}_{\bar k}, \ I\rightarrow \bar I}
\left| \frac{e^{i \bar k \cdot g(I) T} - 1}{i \bar k \cdot g(I)T} \right| = \lim_{J \to
  0} \left| \frac{e^{i J T} - 1}{i J T} \right| = \lim_{J \to 0}
\frac{\sqrt{2[1 - \cos (JT)]}}{|J T|} = 1,
$$
and also 
$$
\sup_{I \in A\backslash {\cal R}_{\bar k}} 
\left| f_{\bar k}(I) \frac{e^{i \bar k \cdot g(I) T}
  - 1}{i \bar k \cdot g(I)T} \right| \ge \lim_{I \notin {\cal R}_{\bar k}, 
\  I\rightarrow \bar I}
\left| f_{\bar k}(I) \frac{e^{i \bar k \cdot g(I) T} - 1}{i \bar k \cdot
  g(I)T} \right| = \norm{ f_{\bar k}(\bar I)} .
$$
As a consequence,
$$
|G^T - \bar{f}|^{\infty} = \sum_{k \in \mathbb{Z}^n} \sup_{I \in A} \left|
(G^T - \bar{f})_k(I) \right| \ge |f_{\bar k}(\bar I)| > 0
$$
that is, $G^T$ does not uniformly Fourier converge to $\bar{f}$
in any set $B\times {\Bbb T}^n$ with $B\subseteq A$ open.\\
\indent We proceed with the same discussion for $F^{\mu}$. By denoting 
$$(F^{\mu} - \bar{f})(I,\varphi) := \sum_{k \in \mathbb{Z}^n} (F^{\mu} - \bar{f})_k(I)e^{ik \cdot \varphi},$$ 
from (\ref{coeff2}) and (\ref{coeff4}) we have
\begin{equation} \label{differenza2}
(F^{\mu} - \bar{f})_k(I) = \begin{cases} 0 & \qquad \text{if } k \cdot g(I) = 0 \\ 
\displaystyle{-\mu \frac{f_k(I)}{i k \cdot g(I) - \mu}} & \qquad \text{if } k \cdot g(I) \ne 0
\end{cases}
\end{equation}
By considering as before $\bar I \in {\cal R}(f)$ such that 
$\bar k \cdot g(\bar I)=0$ and 
$\norm{f_{\bar k}(\bar I)}>0$, for some $\bar k \ne 0$, from
$$
\lim_{I \notin {\cal R}_{\bar k}, 
\ I\rightarrow \bar I} \left| -\mu \frac{f_{\bar k}(I)}{i \bar k \cdot g(I)
  - \mu}\right| = |f_{\bar k}(\bar I)|
$$
we have 
$$
|F^{\mu} - \bar{f}|^{\infty} = \sum_{k \in \mathbb{Z}^n} \sup_{I \in A} \left|
(F^{\mu} - \bar{f})_k(I) \right| \ge |f_{\bar k}(\bar I)| > 0
$$
that is, $F^{\mu}$ does not uniformly Fourier converge to $\bar{f}$
in any set $B\times {\Bbb T}^n$ with $B\subseteq A$ open. \\
\indent We conclude the first part of the proof by showing that also $F^{\mu,\nu}$ does not uniformly Fourier converge to $\bar{f}$. Indeed, in such a case,
formulas (\ref{coeff3}) and (\ref{coeff4}) give
\begin{equation} \label{differenza3}
(F^{\mu,\nu} - \bar{f})_k(I) = 
\begin{cases} \displaystyle{f_k(I) \left[ \frac{\mu}{\mu + \nu \| k \|^2} - 1\right]} & \qquad \text{if } k \cdot g(I) = 0 \\ \\
\displaystyle{-\mu \frac{f_k(I)}{i k \cdot g(I) - \mu - \nu \|k\|^2}} & \qquad \text{if } k \cdot g(I) \ne 0
\end{cases}
\end{equation}
By considering again $\bar k \cdot g(\bar I)=0$ and 
$\norm{f_{\bar k}(\bar I)}>0$ with $\bar k \ne 0$, and  
sequences $(\mu_i,\nu_i) \to 0$, we discuss the following 
two cases.
\begin{enumerate}
\item[$(i)$]  If $\lim_{i \to  +\infty} \frac{\nu_i}{\mu_i} = 0$, we have
$$
\lim_{i \to  +\infty} \lim_{I \notin {\cal R}_{\bar k}, \ I\rightarrow 
\bar I}|(F^{\mu_i,\nu_i} - \bar{f})_{\bar k}(I)| = \lim_{i \to  +\infty} 
\displaystyle{\left|\displaystyle{\mu_i \frac{f_{\bar k}(\bar I)}{\mu_i + \nu_i
      \|\bar k\|}}\right|} = |f_{\bar k}(\bar I)|.
$$
\item[$(ii)$] On the contrary, if the sequence $ \frac{\nu_i}{\mu_i}$ does  
not converge to zero, we consider 
$$
\left|(F^{\mu_i,\nu_i} - \bar{f})_{\bar k}(\bar I)\right| = 
\norm{f_{\bar k}(\bar I)}
 \left| \displaystyle{\left[ \frac{\mu_i}{\mu_i + \nu_i \| \bar k \|^2} - 
1\right]} \right|=\norm{f_{\bar k}(\bar I)} {\nu_i \Norm{\bar k}^2\over
   \mu_i+\nu_i \Norm{\bar k}^2} 
$$
which does not converge to zero as $i$ tends to infinity.    
\end{enumerate} 
As a consequence of all previous cases, we conclude that $F^{\mu, \nu}$ does 
not uniformly Fourier converge to $\bar{f}$ in any set $B\times {\Bbb T}^n$ 
with $B\subseteq A$ open. \\
\indent We proceed by discussing the convergence to $\bar{f}$ in the $\norm{ \
  \cdot \ }^0$ norm. Since $g: A \to \mathbb{R}^n$ is a diffeomorphism, 
the set of all resonances $\mathcal{R} := \bigcup_{k\in {\Bbb Z}^n\backslash 0} \mathcal{R}_k$  
has measure zero. Consequently, the norm $\norm{ \ \cdot \ }^0$ can be 
rewritten as
$$
|u|^0 = \sum_{k \in \mathbb{Z}^n} \int_{\tilde{A}} |u_k(I)| dI
$$
where
\begin{equation} \label{A TILDE}
\tilde{A} := A\backslash {\cal R}= \{I \in A: \ k \cdot g(I) \ne 0 \text{ for
  all } k \in \mathbb{Z}^n\backslash 0\}.
\end{equation}
We first prove $\lim_{T \to +\infty} |G^T - \bar{f}|^0 = 0$. From (\ref{coeff1}) and (\ref{coeff4}) we immediately obtain 
$$
(G^T - \bar{f})_k(I) = f_k(I) \frac{e^{ik \cdot g(I)} - 1}{ik \cdot g(I)T}
\qquad \forall I \in \tilde{A}
$$
so that
$$
\int_{\tilde{A}} |(G^T - \bar{f})_k(I)| dI = \int_{\tilde{A}} \left| f_k(I)
\frac{e^{i k \cdot g(I) T} - 1}{i k \cdot g(I) T} \right| dI
$$
$$
\le |f_k|^{\infty} \int_{\tilde{A}} \frac{\sqrt{\sin^2(k \cdot g(I) T) +
    [\cos(k \cdot g(I) T) - 1]^2}}{|k \cdot g(I)|T} dI = |f_k|^{\infty}
\int_{\tilde{A}} \frac{\sqrt{2[1-\cos(k \cdot g(I) T)]}}{|k \cdot g(I)|T} dI.
$$
Using the change of variables 
\begin{equation} \label{first change}
I \mapsto J := g(I) 
\end{equation}
and assumption (\ref{GI}), we obtain
\begin{equation}
\int_{\tilde{A}} |(G^T - \bar{f})_k(I)| dI  \le |f_k|^{\infty}
\int_{\tilde{A}} \frac{\sqrt{2[1-\cos(k \cdot g(I) T)]}}{|k \cdot g(I)|T} dI
\le \frac{|f_k|^{\infty}}{m} \int_{g(\tilde A)} \frac{\sqrt{2[1 - \cos (k
      \cdot J T)]}}{|k \cdot J|T} dJ .
\label{intjk}
\end{equation}
Let now $\tilde{e}_1, \ldots , \tilde{e}_n$ be an orthonormal basis of
$\mathbb{R}^n$ with $k \in \langle \tilde{e}_2, \ldots , \tilde{e}_n
\rangle^{\bot}$ and $R$ a rotation matrix such that $Rk = \|k\| \tilde{e}_1$
(the dependence of the basis and the rotation matrix on $k \in \mathbb{Z}^n$
is here omitted). By the further change of variables
\begin{equation} \label{second change}
J \mapsto x := RJ  
\end{equation}
the quantity $k\cdot J$ in (\ref{intjk}) becomes $k \cdot J =\Norm{k}x_1$, 
and for any $x$ in the integration domain $Rg(\tilde A)$ we have 
$x= Rg(I)$ with $I\in \tilde A$ and $\Norm{x}\leq \Norm{g(I)} \leq C$. As a consequence, we obtain
$$
\int_{\tilde{A}} |(G^T - \bar{f})_k(I)| dI  \le  
\frac{|f_k|^{\infty}}{m} \int_{g(\tilde{A})} 
\frac{\sqrt{2[1 - \cos (k \cdot J T)]}}{|k \cdot J|T} dJ =
\frac{|f_k|^{\infty}}{m} \int_{Rg(\tilde{A})} \frac{\sqrt{2[1 - \cos (\|k\|
      x_1 T)]}}{\|k\||x_1|T} dx_1 \ldots dx_n
$$ 
$$
\le \frac{|f_k|^{\infty} C^{n-1}}{m} \int^{C}_{-C} \frac{\sqrt{2[1 - \cos
      (\|k\| x_1 T)]}}{\|k\||x_1|T} dx_1 =  \frac{|f_k|^{\infty} C^{n-1}}{m \|
  k \| T} \int^{\|k\| C T}_{-\| k \| C T} \frac{\sqrt{2(1 - \cos y)}}{|y|} dy
$$
$$
= \frac{2 |f_k|^{\infty} C^{n-1}}{m \| k \| T} \int^{\|k\| C T/2}_{-\| k \| C T/2} \left| \frac{ \sin y}{y} \right| dy = \frac{4 |f_k|^{\infty} C^{n-1}}{m \| k \| T} \int^{\|k\| C T/2}_{0} \norm{\frac{ \sin y}{y}} dy$$
\begin{equation} \label{prima-stima}
\le \frac{4 |f_k|^{\infty} C^{n-1}}{m \| k \| T} \int^{2\pi}_0 
\norm{\frac{\sin y}{y}} dy + \frac{4 |f_k|^{\infty} C^{n-1}}{m \| k \| T}
\int_{2\pi}^{\|k\|TC/2} \frac{1}{y} dy \leq \frac{4 |f_k|^{\infty} C^{n-1}}{m \|
  k \| T} [l_0 + \log (\| k\| T C)]
\end{equation}
with $l_0:=\int^{2\pi}_0 \norm{\frac{\sin y}{y}} dy\leq 3$. Consequently,
\begin{equation*} 
|G^T - \bar{f}|^0 \le \sum_{k\in {\Bbb Z}^n\backslash 0} \frac{4 |f_k|^{\infty} C^{n-1}}{m \| k \| T} [3 + \log(\|k\| T C)]
\end{equation*}
proving that $G^T$ converges to $\bar{f}$ in the $\norm{\ \cdot \ }^0$ norm. \\
\indent We conclude the proof with the convergence of $F^{\mu}$ and $F^{\mu,\nu}$ to $\bar{f}$. By using formulas (\ref{coeff3}) and (\ref{coeff4}), we have
$$(F^{\mu,\nu} - \bar{f})_k(I) = - \mu \frac{f_k(I)}{i k \cdot g(I) - \mu - \nu \| k \|^2} \qquad \forall I \in \tilde{A}.$$
Hence
$$\int_{\tilde{A}}|(F^{\mu,\nu} - \bar{f})_k(I)| dI = \mu \int_{\tilde{A}} \frac{|f_k(I)|}{\sqrt{(\mu + \nu \|k\|^2)^2 + (k \cdot g(I))^2}} dI$$
$$\le \mu |f_k|^{\infty} \int_{\tilde{A}} \frac{1}{\sqrt{(\mu + \nu \|k\|^2)^2 + (k \cdot g(I))^2}} dI$$
The same changes of variables of the previous case, see (\ref{first change}) and (\ref{second change}), provide
$$\int_{\tilde{A}}|(F^{\mu,\nu} - \bar{f})_k(I)| dI \le \frac{\mu |f_k|^{\infty} C^{n-1}}{m} \int_{-C}^{C} \frac{1}{\sqrt{(\mu + \nu \|k\|^2)^2 + \|k\|^2x_1^2}} dx_1$$ 
$$= \frac{\mu |f_k|^{\infty} C^{n-1}}{m\|k\|} \int^{\frac{\|k\| C}{\mu + \nu \|k\|^2}}_{-\frac{\|k\| C}{\mu + \nu \|k\|^2}} \displaystyle{\frac{1}{\sqrt{1 + x^2}}dx} =  \frac{2\mu |f_k|^{\infty} C^{n-1}}{m\|k\|} \int^{\frac{\|k\| C}{\mu + \nu \|k\|^2}}_{0} \displaystyle{\frac{1}{\sqrt{1 + x^2}}dx}$$
$$= \frac{2\mu |f_k|^{\infty} C^{n-1}}{m\|k\|} \left[ \int^{1}_0 \frac{1}{\sqrt{1 + x^2}} dx + \int^{\frac{\|k\| C}{\mu + \nu \|k\|^2}}_{1} \displaystyle{\frac{1}{\sqrt{1 + x^2}}dx} \right] \le \frac{2\mu |f_k|^{\infty} C^{n-1}}{m\|k\|} \left[ l_1 + \int^{\frac{\|k\| C}{\mu + \nu \|k\|^2}}_{1} \displaystyle{\frac{1}{x}dx} \right]$$
\begin{equation} \label{seconda stima}
= \frac{2\mu |f_k|^{\infty} C^{n-1}}{m\|k\|} \left[ l_1 + 
\log {\|k\| C \over \mu + \nu \|k\|^2} \right]
\end{equation}
with $l_1:= {\rm arcsinh}\ 1 \leq 1 $. Consequently
$$
|F^{\mu,\nu} - \bar{f}|^0 \le {2 C^{n-1}\over m}
\sum_{k\in {\Bbb Z}^n \backslash 0}|f_k|^{\infty}{\mu \over \Norm{k}}
\left[ 1 + \log {\|k\| C \over \mu + \nu \|k\|^2} \right]
$$
and for $\nu = 0$ we obtain also (\ref{convfm}).  

Inequalities  (\ref{convfmn}) and (\ref{convfm})
respectively prove that $F^{\mu,\nu}$ converges 
to $\bar{f}$ for $(\mu, \nu) \to (0,0)$ and $F^{\mu}$ converges to $\bar{f}$ for $\mu \to 0$ in the $\norm{ \ \cdot \ }^0$ norm. \hfill $\Box$

\subsection{Proof of Proposition \ref{PROP1}}

Let us consider any open set $B\subseteq A$. 
Since $g: A \to \mathbb{R}^n$ is a  diffeomorphism, the 
$\norm{ \ \cdot \ }^1$ norm in $B\times {\Bbb T}^n$ 
--see (\ref{norm3-intro})-- can be rewritten as
$$
|u|^1= \sum_{k \in \mathbb{Z}^n} \left\{ \int_{\tilde{B}} |u_k(I)| dI +
\sum_{j=1}^n \int_{\tilde{B}} \left( \norm{ \frac{\partial u_k}{\partial
  I_j}(I)} + \norm{k_j u_k(I)} \right) dI \right\}
$$
with $\tilde{B}= B \cap {\tilde A}$, see (\ref{A TILDE}). \\
\indent We first prove that $G^T$ does not converge to $\bar{f}$ in the 
set $B\times {\Bbb T}^n$. It is sufficient 
to prove that there exists $\epsilon>0$
such that for any large $T$ we have
\begin{equation}
\sum_{j=1}^n\sum_{k \in \mathbb{Z}^n} \int_{\tilde{B}} \left| \left(\frac{\partial G_k^T}{\partial I_j} - \frac{\partial \bar{f}_k}{\partial I_j}\right)(I)  \right| dI
 > \epsilon  .
\label{lowb}
\end{equation}
From  (\ref{coeff1}) and (\ref{coeff4}), for 
any $I \in \tilde{B}$ we have
$$
(G^T-\bar{f})_k(I) = 
\begin{cases}  
\displaystyle{f_k(I) \frac{e^{ik \cdot g(I) T} - 1}{ik \cdot g(I) T}} & \qquad \text{if } 0 \ne k \in \mathbb{Z}^n \\ 
0 & \qquad \text{if } k = 0
\end{cases}
$$
so that
\begin{equation} \label{diffe}
\left(\frac{\partial G_k^T}{\partial I_j} - \frac{\partial \bar{f}_k}{\partial I_j}\right)(I) = 
\begin{cases} 
\displaystyle{\frac{\partial f_k}{\partial I_j}(I) \frac{e^{i k \cdot g(I) T} - 1}{i k \cdot g(I)T} + f_k(I) \frac{\partial}{\partial I_j} \left(\frac{e^{i k \cdot g(I) T} - 1}{i k \cdot g(I)T}\right)} & \qquad \text{if } 0 \ne k \in \mathbb{Z}^n \\ 
0 & \qquad \text{if } k = 0
\end{cases} \end{equation}
We notice that the first addendum in (\ref{diffe}) tends to $0$, that is 
$$\lim_{T \to +\infty} \int_{\tilde{B}} \left| \frac{\partial f_k}{\partial I_j}(I) \frac{e^{i k \cdot g(I) T} - 1}{i k \cdot g(I)T} \right| dI = 0.$$ 
Indeed, by using the changes of variables (\ref{first change}) and (\ref{second change}) as in the proof of Proposition \ref{PROP2} --see also (\ref{prima-stima})-- we obtain 
$$
\int_{\tilde{B}} \left| \frac{\partial f_k}{\partial I_j}(I) \frac{e^{i k
    \cdot g(I) T} - 1}{i k \cdot g(I)T} \right| dI \le 
 \left|\frac{\partial f_k}{\partial I_j}\right|^{\infty}
\frac{4 C^{n-1}}{m \| k \| T} [l_0 + \log (\| k\| T C)].
$$
As a consequence, it remains to study the other term of the equality
(\ref{diffe}), precisely
$$
\sum_{j=1}^n \sum_{k\in {\Bbb Z}^n} \int_{\tilde{B}} 
\norm{f_k(I)} \norm{\frac{\partial}{\partial I_j} \left(\frac{e^{i k \cdot
      g(I) T} - 1}{i k \cdot g(I)T}\right)}dI  
$$
$$
=\sum_{j=1}^n \sum_{k\in {\Bbb Z}^n}  \int_{\tilde{B}} 
\norm{f_k(I)}\norm{ \frac{\partial}{\partial I_j}(i k \cdot g(I))}
{\sqrt{2+(k \cdot g)^2 T^2 -2 k \cdot g T \sin(k \cdot g T)
-2\cos(k \cdot g T)}\over (k\cdot g)^2 T}dI  
$$
$$
= \sum_{k\in {\Bbb Z}^n} \int_{\tilde{B}}  
\norm{f_k(I)}\Norm{  {\partial g \over \partial I}^T   k}_1
{\sqrt{2+(k \cdot g)^2 T^2 -2 k \cdot g T \sin(k \cdot g T)
-2\cos(k \cdot g T)}\over (k\cdot g)^2 T}dI 
$$
where
$$
\Norm{  {\partial g \over \partial I}^T   k}_1:= 
\sum_{j=1}^n \norm{ \Big ( {\partial g \over \partial I}^T   k 
\Big )_j}=
\sum_{j=1}^n \norm{\sum_{i=1}^n  {\partial g_i \over \partial I_j} k_i}=
\sum_{j=1}^n \norm{ \frac{\partial}{\partial I_j}(i k \cdot g(I))}   .
$$
Since ${\cal R}(f)$ is dense in $A$, there exists $\bar I\in  B \cap 
{\cal R}(f)$ such that $k \cdot g(\bar I)=0$ and $\norm{f_k(\bar I)}>0$ 
for some $k \in {\Bbb Z}^n\backslash 0$. In particular, there exist 
$\delta, \lambda_1,\lambda_2 >0$ 
 (independent of $T$) such that 
the closed ball
$$
B_\delta(\bar I) = \{ I: \Norm{I-\bar I}\leq \delta\} 
$$ 
is contained in $B$, and also for any $I\in B_\delta(\bar I)$ we have
\begin{equation*}
\norm{f_k(I)} \geq \lambda_1
\end{equation*}
and 
\begin{equation*}
\min_{ \Norm{u}=1} \Norm{ {\partial g \over \partial I}^T   u}_1 
 \geq \lambda_2  .
\end{equation*}
Let us remark that the constant $\lambda_1$ satisfies 
$0 <  \lambda_1 \leq \norm{f_k}^\infty$. 
The constant $\lambda_2$ is indeed strictly positive, since otherwise there 
would exist $u \ne 0$ with ${\partial g \over \partial I}^T   u=0$, which  
is in contradiction with (\ref{GI}). From (\ref{GI}), there exists also a 
constant $M>0$ such that 
\begin{equation}
\norm{\det {\partial g\over \partial I}(I)} \leq M  
\label{detgi}
\end{equation}
for any $I\in A$. As a consequence, we have
$$
\sum_{j=1}^n\sum_{{\tilde k}\in {\Bbb Z}^n} \int_{\tilde{B}} 
\norm{f_{\tilde k}(I)} \norm{\frac{\partial}{\partial I_j} 
\left(\frac{e^{i {\tilde k} \cdot
      g(I) T} - 1}{i {\tilde k} \cdot g(I)T}\right)}dI  
$$
$$
\geq \lambda_1 \lambda_2 \Norm{k} 
 \int_{ B_\delta(\bar I)} 
{\sqrt{2+(k \cdot g)^2 T^2 -2 k \cdot g T \sin(k \cdot g T)
-2\cos(k \cdot g T)}\over (k\cdot g)^2 T}dI  .
$$
By performing the change of variables $J:=g(I)$ and using (\ref{detgi}), 
the above term has the lower bound 
$$
{\lambda_1\lambda_2\over M}\Norm{k}
\int_{ g( B_\delta(\bar I))}  
{\sqrt{2+(k \cdot J)^2 T^2 -2 k \cdot J T \sin(k \cdot J T)
-2\cos(k \cdot J T)}\over (k\cdot J)^2 T}dJ  ,
$$
which, using the additional change of variables $x:= R J$ as 
in  (\ref{second change}), equals to
$$
{\lambda_1\lambda_2\over M}\Norm{k}
\int_{R g( B_\delta(\bar I) )} 
{\sqrt{2+\Norm{k}^2 x_1^2 T^2 -2 \Norm{k}x_1 T \sin(\Norm{k}x_1  T)
-2\cos(\Norm{k}x_1  T)}\over \Norm{k}^2x_1^2 T}dx  .
$$
We consider $\tilde \delta>0$ possibly depending 
on $k,\bar I$ (but independent of $T$) such that 
$$
\left \{ x: \max_{j=1,\ldots ,n} \norm{x_j - Rg(\bar I)_j} \leq \tilde \delta
\right \} \subseteq Rg(B_\delta (\bar I))  ,
$$
so that we have 
$$
{\lambda_1\lambda_2\over M}\Norm{k}
\int_{R g( B_\delta(\bar I) )} 
{\sqrt{2+\Norm{k}^2 x_1^2 T^2 -2 \Norm{k}x_1 T \sin(\Norm{k}x_1  T)
-2\cos(\Norm{k}x_1  T)}\over \Norm{k}^2x_1^2 T}dx 
$$
$$
\geq 
{\lambda_1\lambda_2\over M}\Norm{k} {\tilde \delta}^{n-1}
\int_{Rg(\bar I)_1 - \tilde \delta}^{Rg(\bar I)_1 + \tilde \delta}
{\sqrt{2+\Norm{k}^2 x_1^2 T^2 -2 \Norm{k}x_1 T \sin(\Norm{k}x_1  T)
-2\cos(\Norm{k}x_1  T)}\over \Norm{k}^2x_1^2 T}dx_1
$$
$$
={\lambda_1\lambda_2\over M} {\tilde \delta}^{n-1}
\int_{\Norm{k}T(  Rg(\bar I)_1 - \tilde \delta)}^{
\Norm{k}T(Rg(\bar I)_1 + \tilde \delta)}
{\sqrt{2+y^2 -2 y  \sin y-2\cos y}\over y^2}dy  .
$$
We remark that, since the change of variables (\ref{second change}) is
performed by a matrix $R$ such that $R k =\Norm{k}\tilde e_1$, so that
$$
Rg(\bar I)_1 = \tilde e_1 \cdot Rg(\bar I) = {1\over \Norm{k}}
Rk \cdot R g(\bar I) = {1\over \Norm{k}} k\cdot g(\bar I)=0  ,
$$
we have
$$
{\lambda_1\lambda_2\over M} {\tilde \delta}^{n-1}
\int_{\Norm{k}T(  Rg(\bar I)_1 - \tilde \delta)}^{
\Norm{k}T(Rg(\bar I)_1 + \tilde \delta)}
{\sqrt{2+y^2 -2 y  \sin y-2\cos y}\over y^2}dy = 
{\lambda_1\lambda_2\over M} {\tilde \delta}^{n-1}
\int_{-\Norm{k}T \tilde \delta}^{
\Norm{k}T\tilde \delta}
{\sqrt{2+y^2 -2 y  \sin y-2\cos y}\over y^2}dy  .
$$
Since for any $y\in {\Bbb R}$ we have
$$
2+y^2 -2 y  \sin y-2\cos y \geq {y^4 \over 4(1+y^2)}  ,
$$
we conclude 
$$
{\lambda_1\lambda_2\over M} {\tilde \delta}^{n-1}
\int_{-\Norm{k}T\tilde \delta}^{
\Norm{k}T\tilde \delta}
{\sqrt{2+y^2 -2 y  \sin y-2\cos y}\over y^2}dy \geq 
{\lambda_1\lambda_2\over 2 M} {\tilde \delta}^{n-1}
\int_{-\Norm{k}T\tilde \delta}^{
\Norm{k}T\tilde \delta}{1\over \sqrt{1+y^2}}dy 
$$
\begin{equation*}
= {\lambda_1\lambda_2\over M} {\tilde \delta}^{n-1}
\int_{0}^{
\Norm{k}T\tilde \delta}{1\over \sqrt{1+y^2}}dy 
={\lambda_1\lambda_2\over M} {\tilde \delta}^{n-1} {\rm  arcsinh}
(\Norm{k}T\tilde \delta)  .
\end{equation*}
Since
$$
\lim_{T\rightarrow +\infty} {\rm  arcsinh}
(\Norm{k}T\tilde \delta)= + \infty  ,
$$
with a suitable definition of $\epsilon$, one immediately obtains  
(\ref{lowb}).\\
\indent We proceed by proving that $F^{\mu}$ does not converge to $\bar{f}$
in the set $B\times {\Bbb T}^n$. It is sufficient 
to prove that there exists $\epsilon>0$ such that for any small $\mu$ we have
\begin{equation}
\sum_{j=1}^n\sum_{k \in \mathbb{Z}^n} \int_{\tilde{B}}
 \left| \left(\frac{\partial F_k^{\mu}}{\partial I_j} -
 \frac{\partial \bar{f}_k}{\partial I_j}\right)(I)  \right| dI > \epsilon .
\label{lowbi}
\end{equation}
From (\ref{differenza2}), for any $I \in \tilde{B}$ we have
$$
(F^{\mu} - \bar{f})_k(I) = 
\begin{cases} 
F_k^{\mu}(I) = \displaystyle{-\frac{\mu f_k(I)}{i k \cdot g(I) - \mu}} & \qquad \text{if } 0 \ne k \in \mathbb{Z}^n \\
0  & \qquad \text{if } k = 0
\end{cases}
$$
so that we have to estimate 
$$
\sum_{j=1}^n\sum_{k \in \mathbb{Z}^n\backslash 0} \int_{\tilde{B}}
 \left| \frac{\partial F_k^{\mu}}{\partial I_j}(I)  \right| dI  .
$$
By direct computations we obtain
\begin{eqnarray}
\sum_{j=1}^n\sum_{k \in \mathbb{Z}^n\backslash 0}\int_{\tilde{B}} \left| \frac{\partial F^{\mu}_k}{\partial I_j}(I)\right| dI 
&=&\sum_{j=1}^n\sum_{k \in \mathbb{Z}^n\backslash 0} \mu \int_{\tilde{B}} 
\frac{1}{\norm{i k \cdot g(I) - \mu}^2} 
\left| \frac{\partial f_k}{\partial I_j}(I) (i k \cdot g(I) - \mu) - f_k(I)
\frac{\partial}{\partial I_j}(i k \cdot g(I)) \right|
dI\cr
&=&\sum_{j=1}^n\sum_{k \in \mathbb{Z}^n\backslash 0} \mu \int_{\tilde{B}}  {1\over \norm{i k \cdot g(I) - \mu}^2}
\sqrt{\mu^2 \Big ( \frac{\partial f_k}{\partial I_j}\Big )^2 +
\Big (  k \cdot g(I) \frac{\partial f_k}{\partial I_j}- f_k 
k \cdot {\partial g\over \partial I_j}\Big )^2}dI\cr
&\geq& \sum_{j=1}^n\sum_{k \in \mathbb{Z}^n\backslash 0}  \mu \int_{\tilde{B}}  {1\over \norm{i k \cdot g(I) - \mu}^2}
\norm{  k \cdot g(I) \frac{\partial f_k}{\partial I_j}- f_k 
k \cdot {\partial g\over \partial I_j}}dI  .
\label{sumint}
\end{eqnarray}
As before, we consider $\bar I \in \tilde{B}\cap {\cal R}(f)$, so that 
there exists  $k \in {\Bbb Z}^n$ such that $k \cdot g(\bar I)=0$ and 
$\norm{f_k(\bar I)}>0$. In particular, there exist 
$\delta, \lambda_1,\lambda_2 >0$ (independent of $T$) such that 
the closed ball
$$
B_\delta(\bar I) = \{ I: \Norm{I-\bar I}\leq \delta\} 
$$ 
is contained in $B$, and also for any $I\in B_\delta(\bar I)$ we have
\begin{equation*}
\norm{f_k(I)} \geq \lambda_1  ,
\end{equation*}
and 
\begin{equation*}
\min_{ \Norm{u}=1} \Norm{ {\partial g \over \partial I}^T   u}_1 
 \geq \lambda_2  .
\end{equation*}
Since $\lambda>0$ is a Lipschitz constant for $g$ in the set $A$, 
for any $I\in B_\delta(\bar I)$ we also have 
\begin{equation*}
\norm{k \cdot g(I)} \leq \Norm{k} \lambda \delta  .
\end{equation*} 
The series in (\ref{sumint}) has therefore the lower bound
$$
\mu \sum_{j=1}^n\sum_{{\tilde k} \in \mathbb{Z}^n\backslash 0}   
\int_{\tilde{B}}  
{1\over \norm{i {\tilde k} \cdot g(I) - \mu}^2}
\norm{  {\tilde k} \cdot g(I) \frac{\partial f_{\tilde k}}{\partial I_j}- 
f_{\tilde k} 
{\tilde k} \cdot {\partial g\over \partial I_j}}dI 
$$
$$
\geq \mu \sum_{j=1}^n  \int_{B_\delta(\bar I)}
{1\over \norm{i k \cdot g(I) - \mu}^2}
\norm{  k \cdot g(I) \frac{\partial f_k}{\partial I_j}- f_k 
k \cdot {\partial g\over \partial I_j}}dI
$$
$$
\geq   \mu \int_{B_\delta(\bar I) }  {1\over \norm{i k \cdot g(I) - \mu}^2} 
\norm{f_k} \Norm{  {\partial g \over \partial I}^T   k}_1dI -
\mu\sum_{j=1}^n  \int_{B_\delta(\bar I) }  
{1\over \norm{i k \cdot g(I) - \mu}^2}  
\norm{k\cdot g}\norm{\frac{\partial f_k}{\partial I_j}}dI
$$
$$
\geq \lambda_1\lambda_2 \Norm{k} \mu 
\int_{B_\delta(\bar I) }  
{1\over \norm{i k \cdot g(I) - \mu}^2}dI - \mu \Norm{k}\delta \lambda
\Big (\sum_{j=1}^n \norm{\frac{\partial f_k}{\partial I_j}}^\infty \Big )
\int_{B_\delta(\bar I) } {1\over \norm{i k \cdot g(I) - \mu}^2}dI  
$$
$$
= \left (\lambda_1\lambda_2 - \delta \lambda 
\sum_{j=1}^n \norm{\frac{\partial f_k}{\partial I_j}}^\infty \right )
\Norm{k} \mu  \int_{B_\delta(\bar I) } {1\over \norm{i k \cdot g(I) -
    \mu}^2}dI  .
$$
First, we remark that in the case $\sum_{j=1}^n \norm{\frac{\partial
    f_k}{\partial I_j}}^\infty>0$,  
it is not restrictive to choose $\delta$ satisfying
$$
\delta \leq {\lambda_1\lambda_2 \over 2 \lambda 
\sum_{j=1}^n \norm{\frac{\partial f_k}{\partial I_j}}^\infty}  ,
$$
so that we have
$$
\sum_{j=1}^n\sum_{{\tilde k} \in \mathbb{Z}^n\backslash 0}  \mu \int_{\tilde{B}}  
{1\over \norm{i {\tilde k} \cdot g(I) - \mu}^2}
\norm{  {\tilde k} \cdot g(I) \frac{\partial f_{\tilde k}}{\partial I_j}- 
f_{\tilde k}
{\tilde k} \cdot {\partial g\over \partial I_j}}dI 
\geq {\lambda_1\lambda_2 \over 2} \Norm{k} \mu  \int_{B_\delta(\bar I) } {1\over \norm{i k \cdot g(I) -
    \mu}^2}dI  .
$$
Then by performing the change of variables $J:=g(I)$ and using (\ref{detgi}) 
we obtain the lower bound 
$$
 {\lambda_1\lambda_2 \over 2} \Norm{k} \mu  \int_{B_\delta(\bar I) } {1\over \norm{i k \cdot g(I) -
    \mu}^2}dI \geq 
 {\lambda_1\lambda_2 \over 2 M} \Norm{k} \mu 
 \int_{g(B_\delta(\bar I)) }  {1\over \norm{i k \cdot J -
    \mu}^2}dJ
$$
$$
= {\lambda_1\lambda_2 \over 2 M} \Norm{k} \mu 
 \int_{g(B_\delta(\bar I)) }  {1\over  \sqrt{(k\cdot J)^2 +\mu^2}}dJ
$$
which, by the additional change of variables $x:= R J$ as 
in  (\ref{second change}), can be written as
$$
 {\lambda_1\lambda_2 \over 2 M} \Norm{k} \mu 
 \int_{R g(B_\delta(\bar I)) }   {1\over  \sqrt{\Norm{k}^2x_1^2 +\mu^2}}
dx  .
$$
Since there exists $\tilde \delta>0$ possibly depending 
on $k,\bar I$ (but independent of $\mu$) such that 
$$
\left \{ x: \max_{j=1,\ldots ,n} \norm{x_j - Rg(\bar I)_j} \leq \tilde \delta
\right \} \subseteq Rg (B_\delta (\bar I))  ,
$$
we obtain the lower bound
$$
 {\lambda_1\lambda_2 \over 2 M} \Norm{k} \mu 
 \int_{R g(B_\delta(\bar I)) }   {1\over  \sqrt{\Norm{k}^2x_1^2 +\mu^2}}dx 
\geq  {\lambda_1\lambda_2 \over 2 M} \Norm{k} {\tilde \delta}^{n-1}
\mu  \int_{- \tilde \delta}^{ \tilde \delta}
 {1\over  \sqrt{\Norm{k}^2x_1^2 +\mu^2}}dx_1
$$
$$
= {\lambda_1\lambda_2 \over 2 M}
\int_{-{\Norm{k}\over \mu}\tilde \delta}^{
{\Norm{k}\over \mu} \tilde \delta}
{1\over 1+ y^2}dy  = {\lambda_1\lambda_2 \over M} 
\arctan {\Norm{k}\over \mu}\tilde \delta  .
$$
Since we have
$$
\lim_{\mu \rightarrow 0^+} \arctan {\Norm{k}\over \mu}\tilde \delta ={\pi
  \over 2}  ,
$$
with a suitable definition of $\epsilon$  one immediately obtains  
(\ref{lowbi}).

\indent We conclude our proof by showing the convergence of $F^{\mu,\nu}$ to
$\bar{f}$ in $A\times {\Bbb T}^n$ on sequences $(\mu_i,\nu_i) \to (0,0)$ such 
that
\begin{equation}
\lim_{i \to 0} \frac{\mu_i}{\nu_i} = 0 .
\end{equation}
We first provide an estimate on the different contributions to
$$
\norm{F^{\mu,\nu}-\bar f}^1 = 
\norm{F^{\mu,\nu}-\bar f}^0+\sum_{k\in {\Bbb Z}^n}
\sum_{j=1}^n \int_A \norm{{\partial\over \partial I_j}
(F^{\mu,\nu} - \bar{f})_k}dI
+\sum_{k\in {\Bbb Z}^n}
\sum_{j=1}^n \int_A \norm{k_j}\norm{(F^{\mu,\nu} - \bar{f})_k}dI  .
$$
The first term $\norm{F^{\mu,\nu}-\bar f}^0$ has been already 
estimated (see (\ref{convfmn}))
\begin{equation}
 |F^{\mu,\nu} - \bar{f}|^0 \le  {2  C^{n-1} \over m} 
\sum_{k\in {\Bbb Z}^n\backslash 0 }|f_k|^{\infty}\ \frac{\mu }{\|k\|} \left( 1 + \log {\|k\| C \over \mu+\nu\Norm{k}^2 }\right ) .
\label{c1}
\end{equation}
Then, for any $I\in \tilde{A}$, from 
(\ref{differenza3}) we have 
$$
(F^{\mu,\nu} - \bar{f})_k(I) = \begin{cases}
F_k^{\mu,\nu}(I) = \displaystyle{-\mu \frac{f_k(I)}{i k \cdot g(I) - \mu - \nu \|k\|^2}}  & \qquad \text{if } 0 \ne k \in \mathbb{Z}^n \\
0  & \qquad \text{if } k = 0
\end{cases}
$$
so that we need to estimate 
$$
\int_{\tilde{A}} \left| \frac{\partial F^{\mu,\nu}_k}{\partial I_j}(I) \right|
dI =  \mu \int_{\tilde{A}} \left| \frac{\partial f_k}{\partial I_j}(I)
\frac{1}{i k \cdot g(I) - \mu -\nu \|k\|^2} + f_k(I) \frac{\partial}{\partial
  I_j} \left( \frac{1}{i k \cdot g(I) - \mu - \nu \|k\|^2} \right) \right| dI
$$
for any $k\in {\Bbb Z}^n\backslash 0$. By using the changes of variables (\ref{first change}) and
(\ref{second change}) as in the proof of Proposition \ref{PROP2} --and  
proceeding as in estimate (\ref{seconda stima})-- we obtain 
\begin{equation} 
\mu \sum_{k\in {\Bbb Z}^n \backslash 0}\sum_{j=1}^n 
\int_{\tilde{A}} \left| \frac{\partial
  f_k}{\partial I_j}(I) \frac{1}{i k \cdot g(I) - \mu -\nu_i \|k\|^2} \right| dI
\le {2  C^{n-1}\over m}\sum_{k\in {\Bbb Z}^n \backslash 0} 
\left (\sum_{j=1}^n  \left|\frac{\partial
    f_k}{\partial I_j}\right|^{\infty}\right ) {\mu \over \Norm{k}}
\left[ 1 + \log  {\|k\| C \over  \mu + \nu \|k\|^2}\right ]  .
\label{c2}
\end{equation}
Using (\ref{GI}) we first obtain
$$
 \sum_{k\in {\Bbb Z}^n \backslash 0}\sum_{j=1}^n \mu \int_{\tilde{A}} \left| f_k(I)
\frac{\partial}{\partial I_j} \left( \frac{1}{i k \cdot g(I) - \mu - \nu 
\|k\|^2} \right) \right| dI
$$
$$
\le  \sum_{k\in {\Bbb Z}^n \backslash 0}\sum_{j=1}^n \mu |f_k|^{\infty} \int_{\tilde{A}} \frac{1}{
\norm{i k \cdot g(I) - \mu - \nu \|k\|^2}^2} \left| \frac{\partial}{\partial
  I_j}(i k \cdot g(I)) \right| dI  
$$
$$
\le  \sum_{k\in {\Bbb Z}^n \backslash 0} \mu |f_k|^{\infty} n^2 \Norm{k} D 
\int_{\tilde{A}} \frac{1}{(k\cdot g(I))^2 +(\mu+\nu \|k\|^2)^2} dI  ,
$$
then using the change of variables  (\ref{first change}) and
(\ref{second change}) we have
$$
 \sum_{k\in {\Bbb Z}^n \backslash 0} \mu |f_k|^{\infty} n^2 \Norm{k} D 
\int_{\tilde{A}} \frac{1}{(k\cdot g(I))^2 +(\mu+\nu \|k\|^2)^2}  dI
$$
$$
\leq  \sum_{k\in {\Bbb Z}^n \backslash 0} {n^2 D\over m} |f_k|^{\infty} \Norm{k}  
\mu\int_{g(\tilde{A})} \frac{1}{(k\cdot J)^2 +(\mu+\nu \|k\|^2)^2}  dJ
$$
$$
\leq  \sum_{k\in {\Bbb Z}^n \backslash 0} {n^2C^{n-1} D\over m} |f_k|^{\infty} \Norm{k}  
\mu\int_{-C}^C \frac{1}{\Norm{k}^2 x_1^2 +(\mu+\nu \|k\|^2)^2}  dx_1
$$
$$
\leq  \sum_{k\in {\Bbb Z}^n \backslash 0} {2 n^2C^{n-1} D\over m} |f_k|^{\infty} 
{\Norm{k}   \over (\mu+\nu \|k\|^2)^2}
\mu\int_{0}^C \frac{1}{{\Norm{k}^2\over (\mu+\nu \|k\|^2)^2} x_1^2 +1}  dx_1
$$
$$
=  \sum_{k\in {\Bbb Z}^n \backslash 0} {2 n^2 C^{n-1} D\over m} |f_k|^{\infty}  
{\mu \over (\mu+\nu \|k\|^2)} \int_{0}^{{\Norm{k}C\over \mu+\nu \|k\|^2}}
{1\over 1+y^2}dy
$$
$$
=  \sum_{k\in {\Bbb Z}^n \backslash 0} {2 n^2 C^{n-1} D\over m} |f_k|^{\infty}  
{\mu \over (\mu+\nu \|k\|^2)} \arctan \left( \frac{\|k\| C}{ \mu + \nu
  \|k\|^2} \right).  
$$
From the previous inequality, we obtain
\begin{equation}
\sum_{k\in {\Bbb Z}^n \backslash 0}\sum_{j=1}^n
\mu \int_{\tilde{A}} \left| f_k(I)
\frac{\partial}{\partial I_j} \left( \frac{1}{i k \cdot g(I) - \mu - \nu 
\|k\|^2} \right) \right| dI \leq
{n^2 \pi C^{n-1} D\over m} \sum_{k\in {\Bbb Z}^n \backslash 0}
\norm{f_k}^\infty {\mu \over \mu+\nu \Norm{k}^2}   .
\label{c3}
\end{equation}
In order to conclude the estimate of $\norm{F^{\mu,\nu} - \bar{f}}^1$  it 
remains to consider
$$
\sum_{k\in {\Bbb Z}^n \backslash 0} \sum_{j=1}^n\int_{\tilde{A}}
\left| k_j \frac{\mu_i f_k(I)}{ik \cdot g(I) - \mu_i -\nu_i 
\| k \|^2}\right| dI  .$$
This term is estimated by using the changes of variables (\ref{first change}) and (\ref{second change}), so that
$$
\sum_{k\in {\Bbb Z}^n \backslash 0}\sum_{j=1}^n \mu |k_j| \int_{\tilde{A}} \left| \frac{f_k(I)}{ik \cdot g(I) - \mu -\nu \| k
  \|^2}\right| dI \le \sum_{k\in {\Bbb Z}^n \backslash 0}\sum_{j=1}^n \frac{2 \mu |k_j| 
\left| f_k \right|^{\infty}
  C^{n-1}}{m\|k\|} \left[ l_1 + \log {\|k\| C  \over \mu + \nu \|k\|^2} \right]
$$
\begin{equation}
\leq   {2  n  C^{n-1} \over m}
\sum_{k\in {\Bbb Z}^n \backslash 0}  
\left| f_k \right|^{\infty} 
\mu \left[ 1 + \log {\|k\| C  \over \mu + \nu \|k\|^2 } \right]  .
\label{c4}
\end{equation}
By collecting inequalities (\ref{c1}), (\ref{c2}), (\ref{c3}) and (\ref{c4}), 
 we obtain
$$
\norm{F^{\mu,\nu}-\bar f}^1 \leq {2 C^{n-1}\over m} 
\sum_{k\in {\Bbb Z}^n \backslash 0}  \Bigg ( 
\frac{\mu}{\|k\|} \Big [ 1 + \log {\|k\| C \over \mu+\nu\Norm{k}^2 }\Big ]
\norm{f_k}^\infty 
+  {\mu \over \Norm{k}}
\left[ 1 + \log  {\|k\| C \over  \mu + \nu \|k\|^2}\right ]
\left (\sum_{j=1}^n  \left|\frac{\partial
    f_k}{\partial I_j}\right|^{\infty}\right ) 
$$
$$
+ {1\over 2} n^2 \pi D  {\mu \over \mu+\nu \Norm{k}^2} \norm{f_k}^\infty 
+n \mu \left[ 1 + \log  {\|k\| C \over  \mu + \nu \|k\|^2}\right ] 
\norm{f_k}^\infty  \Bigg )
$$
$$
\leq  {2 C^{n-1}\over m} 
\sum_{k\in {\Bbb Z}^n \backslash 0} \Bigg (  
\mu \Big [ 1 + \log {\|k\| C \over \mu+\nu\Norm{k}^2 }\Big ] 
\Big ( (1+n) \norm{f_k}^\infty  + \sum_{j=1}^n  \left|\frac{\partial
    f_k}{\partial I_j}\right|^{\infty} \Big )+
 {1\over 2} n^2 \pi D  {\mu \over \mu+\nu \Norm{k}^2} \norm{f_k}^\infty\Bigg )
$$
so that (\ref{convfmnorm}) is proved. Since for $\mu,\nu >0$ and 
$\Norm{k}\geq 1$, we have 
$$
\mu  \log {\|k\| C \over \mu+\nu\Norm{k}^2 } \leq 
\mu  \log {C \over \nu} \leq 
{\mu \over \nu} \Big(   \nu  \log {C \over \nu}\Big )
$$
and
$$
{\mu \over \mu+\nu \Norm{k}^2} \leq {\mu \over \nu}  ,
$$
from (\ref{convfmnorm}) we obtain
$$
\norm{F^{\mu,\nu}-\bar f}^1 \leq \Big ( {\mu \over \nu} \Big )
 {2 C^{n-1}\over m}  
\sum_{k\in {\Bbb Z}^n \backslash 0} \Bigg (  
\Big (\nu  +  \nu  \log {C \over \nu} \Big )
\Big ( (1+n) \norm{f_k}^\infty  + \sum_{j=1}^n  \left|\frac{\partial
    f_k}{\partial I_j}\right|^{\infty} \Big )+
 {1\over 2} n^2 \pi D   \norm{f_k}^\infty\Bigg )  .
$$
Therefore, for any sequence $\mu_i,\nu_i>0$ converging to zero 
with $\mu_i/\nu_i$ converging to zero, we have 
$$
\lim_{i\rightarrow +\infty} \norm{F^{\mu_i,\nu_i}-\bar f}^1=0  .
$$
The proof of Proposition \ref{PROP1} is concluded. \hfill $\Box$

\end{document}